\documentclass[12pt]{article}
\usepackage{graphicx,amssymb, vmargin, amsmath, multicol}
\title{The hypermetric cone on seven vertices}
\setpapersize{custom}{21cm}{29.7cm}
\setmarginsrb{1cm}{3cm}{1cm}{3cm}{0pt}{0pt}{0pt}{0pt}

\author{Mathieu Dutour and Michel Deza\footnote{Mathieu.Dutour@ens.fr and Michel.Deza@ens.fr}\\
Laboratoire interdisciplinaire de g\'eom\'etrie appliqu\'ee,\footnote{45 rue d'Ulm, 75005 Paris, France} CNRS/ENS, Paris\\
and Institute of Statistical Mathematics, Tokyo}
\begin{document}
\newcommand{\R}{\ensuremath{\mathbb{R}}}
\newcommand{\N}{\ensuremath{\mathbb{N}}}
\newcommand{\Q}{\ensuremath{\mathbb{Q}}}
\newcommand{\C}{\ensuremath{\mathbb{C}}}
\newcommand{\Z}{\ensuremath{\mathbb{Z}}}
\newcommand{\T}{\ensuremath{\mathbb{T}}}
\newtheorem{theorem}{Theorem}[section]
\newtheorem{definition}[theorem]{Definition}
\newtheorem{conjecture}[theorem]{Conjecture}
\newtheorem{corollary}[theorem]{Corollary}
\newtheorem{proposition}[theorem]{Proposition}
\maketitle

\begin{abstract}
The hypermetric cone $HYP_n$ is the set of vectors $(d_{ij})_{1\leq i< j\leq n}$ satisfying the inequalities
\begin{equation*}
\sum_{1\leq i<j\leq n} b_ib_jd_{ij}\leq 0\mbox{~with~}b_i\in\Z\mbox{~and~}\sum_{i=1}^{n}b_i=1\,.
\end{equation*}
A Delaunay polytope of a lattice is called {\em extremal} if the only affine bijective transformations of it into a Delaunay polytope, are the homotheties; there is correspondance between such Delaunay polytopes and extreme rays of $HYP_n$. We show that unique Delaunay polytopes of root lattices $A_1$ and $E_6$ are the only extreme Delaunay polytopes of dimension at most $6$. We describe also the skeletons and adjacency properties of $HYP_7$ and of its dual.
\end{abstract}

\section{Introduction}
A vector $(d_{ij})_{1\leq i< j\leq n}\in \R^N$ with $N={n\choose 2}$ is called an {\em $n$-hypermetric} if it satisfy the following {\em hypermetric inequalities}:
\begin{equation}\label{Definition-of-hypermetrics}
\sum_{1\leq i<j\leq n} b_ib_jd_{ij}\leq 0\mbox{~with~}b=(b_i)\in\Z^n\mbox{~and~}\sum_{i=1}^{n}b_i=1\,\,.
\end{equation}
The set of vectors satisfying (\ref{Definition-of-hypermetrics}) is called the {\em hypermetric cone} and denoted by $HYP_n$.\\
We have the inclusions $CUT_n\subset HYP_n\subset MET_n$, where $MET_n$ denotes the cone of all semimetrics on $n$ points and $CUT_n$ (see section $3$ below and chapter $4$ of \cite{DL}) is the cone of all semimetrics on $n$ points, which are isometrically embeddable into some space $l^m_1$. In fact, the triangle inequality $d_{ij}\leq d_{ik}+d_{jk}$ is the hypermetric inequality with vector $b$ such that $b_i=b_j=1$, $b_k=-1$ and $b_l=0$, otherwise.\\
For $n\leq 4$ all three cones coincide and $HYP_n=CUT_n$ for $n\leq 6$; so the cone $HYP_7$ is the first proper hypermetric cone\footnote{Apropos, $MET_7$ has $46$ orbits of extreme rays and not $41$ as, by a technical mistake, was given in \cite{GrMet} and \cite{DL}}. See \cite{DL} for detailed study of those cones and their numerous applications in Combinatorial Optimization, Analysis and other areas of mathematics. In particular, the hypermetric cone had direct applications in Geometry of Quadratic Forms; see section \ref{geometry-of-numbers}.\\

In fact, $HYP_n$ is a polyhedral cone (see \cite{DGL93}). Lovasz (see \cite{DL} pp.~201-205) gave another proof of it and the bound $max |b_i|\leq n!2^n
{2n \choose n}^{-1}$ for any vector $b=(b_i)$, defining a facet of $HYP_n$.\\
The group of all permutations on $n$ vertices induces a partition of the set of $i$-dimensional faces of $HYP_n$ into orbits. Baranovskii using his method presented in \cite{B70} found in \cite{Ba} the list of all facets of $HYP_7$: $3773$ facets, divided into $14$ orbits. On the other hand, in \cite{DGL92} were found $29$ orbits of extreme rays of $HYP_7$ by classifying the basic simplexes of the Schlafli polytope of the root lattice $E_6$.\\
In section \ref{Computational-method} we show that the $37170$ extreme rays containing in those $29$ orbits are, in fact, the complete list. It also implies that the Schlafli polytope (unique Delaunay polytope of $E_6$) and the segment $\alpha_1$ (the Delaunay polytope of $A_1$) are only extreme Delaunay polytopes of dimension at most six.\\
In section \ref{geometry-of-hypermetric-cone} we give adjacency properties of the skeletons of $HYP_7$ and of its dual.\\
The computations were done using the programs {\em cdd} (see \cite{Fu}) and {\em nauty} (see \cite{MK}).

\section{Hypermetrics and Delaunay polytopes}\label{geometry-of-numbers}
\noindent For more details on the material of this section see chapters $13$--$16$ of \cite{DL}.\\
Let $L\subset \R^k$ be a $k$-dimensional lattice and let $S=S(c,r)$ be a  sphere in $\R^k$ with center $c$ and radius $r$. Then, $S$ is said to be an {\em empty sphere} in $L$ if the following two conditions hold:
\begin{center}
$\Vert v-c\Vert\geq r$ for all $v\in L$ and the set $S\cap L$ has affine rank $k+1$.
\end{center}
Then, the center $c$ of $S$ is called a {\em hole} in \cite{CS}. The polytope $P$ which is defined as the convex hull of the set $P=S\cap L$ is called a {\em Delaunay polytope}, or (in original terms of Voronoi who introduced them in \cite{Vo}) {\em L-polytope}.\\
On every set $A=\{ v_1, \dots, v_m\}$ of vertices of a Delaunay polytope $P$ we can define a distance function $d_{ij}=\Vert v_i-v_j\Vert^2$. The function $d$ turns out to be a metric and, moreover, a hypermetric. It follows from the following formula (see \cite{As} and \cite{DL} p.~195) :
\begin{equation*}
\sum_{i,j\in A} b_ib_jd_{ij}=2(r^2-\Vert \sum_{i\in A}b_iv_i - c\Vert^2)\leq 0\,.
\end{equation*}
On the other hand, Assouad has shown in \cite{As} that {\em every finite} hypermetric space is a square euclidean distance on a generating set of vertices of a Delaunay polytope of a lattice.\\
For example, in dimension two there are two kinds of combinatorial types of Delaunay polytopes: triangle and rectangle. Since $HYP_3=MET_3$, we see that a triangle can be made a set of vertices of a Delaunay polytope if and only if it has obtuse angles.\\
A Delaunay polytope $P$ is said to be {\em extreme} if the only (up to orthogonal transformations and translations) affine bijective transformations $T$ of $\R^k$ for which $T(P)$ is again a Delaunay polytope, are the homotheties. \cite{DGL92} show that the hypermetric on generating subsets of a extreme Delaunay polytope (see above) lie on extreme ray of $HYP_n$ and that a hypermetric, lying on an extreme ray of $HYP_n$, is the square of euclidean distance on generating subset of extreme Delaunay polytope of dimension at most $n-1$.\\
In \cite{DL} p.~228 there is a more complete dictionary translating the properties of Delaunay polytopes into those of the corresponding hypermetrics.\\
Remind that $E_6$ and $E_7$ are {\em root } lattices defined by
\begin{equation*}
E_6=\{x\in \Z^8\mbox{~:~}x_1+x_2=x_3+\dots+x_8=0\},\,\,E_7=\{x\in \Z^8\mbox{~:~}x_1+x_2+x_3+\dots+x_8=0\}
\end{equation*}
The skeleton of unique Delaunay polytope of $E_6$ is $27$-vertex strongly regular graph, called the {\em Schlafli graph}. In fact, the $29$ orbits of extreme rays of $HYP_7$, found in \cite{DGL92}, were three orbits of extreme rays of $CUT_7$ ({\em cuts}) and $26$ ones, corresponding to all sets of seven vertices of Schlafli graph, which are affine bases (over $\Z$) of $E_6$. The root lattice $E_7$ has two Delaunay polytopes: $1$-simplex and $56$-vertex polytope, called {\em Gosset polytope}. In \cite{Du} were found all $374$ orbits of affine bases for the Gosset polytope.

\section{Computing the extreme rays of $HYP_7$}\label{Computational-method}
\noindent We recall some terminology. Let $C$ be a polyhedral cone 
in $\mathbb{R}^{n}$. Given $v \in \mathbb{R}^{n}$, the inequality 
$\sum_{i=1}^nv_ix_i \geq 0$ is said to be {\it valid} for $C$, if it
 holds for all $x \in C$. Then the set 
$\{ x \in C \vert \sum_{i=1}^nv_ix_i = 0 \}$ is called the {\it face 
of $C$, induced by the valid inequality $\sum_{i=1}^nv_ix_i \geq 0$}. 
A face of dimension $\dim(C) - 1$ is called a {\it facet} of $C$;
 a face of dimension $1$ is called an {\it extreme ray} of $C$.

An extreme ray ray $e$ is said to be {\it incident} to a facet $F$ if 
$e\subset F$. A facet $F$ is said to be {\it incident} to an extreme ray $e$
if $e\subset F$.
if an extreme ray $e$ is {\it incident} to a facet $F$ if $e\subset F$.
Two extreme rays of $C$ are said to be {\it adjacent}, if they
span a two-dimensional face of
$C$. Two facets of $C$ are said to be {\it adjacent}, 
if their intersection has dimension $\dim(C) - 2$.

All $14$ orbits $F_m$, $1\leq m\leq 14$, of facets of $HYP_7$, found by Baranovskii,
are represented below by the corresponding vector $b^{m}$ (see (\ref{Definition-of-hypermetrics})):
\begin{equation*}
\begin{array}{lll}
b^1=(1,1,-1,0,0,0,0);&	b^2=(1,1,1,-1,-1,0,0);&	b^3=(1,1,1,1,-1,-2,0),\\
b^4=(2,1,1,-1,-1,-1,0);&	b^5=(1,1,1,1,-1,-1,-1);&	b^6=(2,2,1,-1,-1,-1,-1),\\
b^7=(1,1,1,1,1,-2,-2),&	b^8=(2,1,1,1,-1,-1,-2); &b^9=(3,1,1,-1,-1,-1,-1),\\
b^{10}=(1,1,1,1,1,-1,-3);&	b^{11}=(2,2,1,1,-1,-1,-3),&	b^{12}=(3,1,1,1,-1,-2,-2),\\
b^{13}=(3,2,1,-1,-1,-1,-2),&	b^{14}=(2,1,1,1,1,-2,-3).
\end{array}
\end{equation*}
It gives the total of $3773$ inequalities. The first ten orbits are the orbits of hypermetric facets of the cut cone $CUT_7$; first four of them come as {\em $0$-extension} of facets of the cone $HYP_6$, i.e. the vector has zero components $x_{ij}$ for some $1\leq i\leq 7$ and all $1\leq j<i$, $i<j\leq 7$. The orbits $F_{11}$--$F_{14}$ consist of some $19$-dimensional simplex faces of $CUT_7$, becoming simplex facets in $HYP_7$.\\
The proof (see \cite{B70} and \cite{BR}) was in terms of volume of simplices; his proof implies that for facet of $HYP_7$ holds $|b_i|\leq 3$ (compare with the bound in introduction). In \cite{BR} the repartitioning polytopes (connected implicitly in \cite{Vo} to facets of $HYP_n$) found for facets of $HYP_7$.\\
Because of the large number of facets of $HYP_7$, it is difficult to find extreme rays just by application of existing programs (see \cite{Fu}). So, let us consider in more detail the cut cone $CUT_7$.\\
Call {\em cut cone} and denote by $CUT_n$ the cone generated by all {\em cuts} $\delta_S$ defined by 
\begin{equation*}
(\delta_S)_{ij}=1\mbox{~if~}\vert S\cap \{i,j\}\vert =1\mbox{~and~} (\delta_S)_{ij}=0 \mbox{,~otherwise},
\end{equation*}
where $S$ is any subset of $\{1,\dots, n\}$.
The cone $CUT_n$ has dimension $n \choose 2$ and $2^{n-1}-1$ non-zero cuts as generators of extreme rays. There are $\lfloor \frac{n}{2}\rfloor$ orbits, corresponding to all non-zero values of $min(|S|,n-|S|)$. The skeleton of $CUT_n$ is the complete graph $K_{2^{n-1}-1}$. See part V of \cite{DL} for a survey on facets of $CUT_n$.\\
The $38780$ facets of the cut cone $CUT_7$ partitioned in $36$ orbits. In \cite{G90} was shown that the known list of $36$ orbits was complete (See \cite{DDL} and chapter $30$ of \cite{DL} for details). Of these $36$ orbits ten are hypermetric. We computed the diameter of the skeleton of dual $CUT_7$: it is exactly $3$ (apropos, the diameter of the skeleton of $MET_n$, $n\geq 4$, is $2$, see \cite{DD1} ). So, we have $CUT_n\subset HYP_n$ and cones $CUT_7$, $HYP_7$ have ten common (hypermetric) facets: $F_1$--$F_{10}$.\\
Each of $26$ orbits of non-hypermetric facets of $CUT_7$ consists of simplex cones, i.e. those facets are incident exactly to $20$ cuts or, in other words, adjacent to $20$ other facets. It turns out that the $26$ orbits of non-hypermetric facets of $CUT_7$ correspond exactly to $26$ orbits of non-cut extreme rays of $HYP_7$.\\
In fact, if $d$ is a point of an extreme ray of $HYP_7$, which is not a cut, then it violates one of the non-hypermetric facet inequalities of $CUT_7$.
More precisely, our computation consist of the following steps:
\begin{enumerate}
\item If $d$ belongs to an non-cut extreme ray of $HYP_7$, then $d\notin CUT_7$.
\item So, there is at least one non-hypermetric facet $F$ of $CUT_7$ with $F(d)<0$.
\item Select a facet $F_i$ for each non-hypermetric orbit of $O_i$ with $1\leq i\leq 26$ and define $26$ subcones $C_i$, $1\leq i\leq 26$, by $C_i=\{d\in HYP_7 \mbox{~:~}F_i(d)\leq 0\}$.
\item The initial set of $3773$ hypermetric inequalities is non-redundant, but adding the inequality $F_i(d)\leq 0$ yields a highly redundant set of inequalities. We remove the redundant inequalities using invariant group (of permutations preserving the cone $C_i$) and linear programming (see polyhedral FAQ\footnote{http://www.ifor.math.ethz.ch/\texttt{\~}fukuda/polyfaq/polyfaq.html} in \cite{Fu}).
\item For each of $26$ subcones we found, by computation, a set of $21$ non-redundant facets, i.e. each of the subcones $C_i$ is a simplex. We get $21$ extreme rays for each of the $26$ subcones.
\item We remove the $20$ extreme rays, which are cuts, from each list and get, for each of these subcones, exactly one non-cut extreme ray.
\end{enumerate}
So, we get an upper bound $26$ for the number of non-cuts orbits of extreme rays. But \cite{DGL92} gave, in fact, a lower bound $26$ for this number. So, we get:
\begin{proposition}
The hypermetric cone $HYP_7$ has $37170$ extreme rays, divided into three orbits corresponding to cuts and $26$ orbits corresponding to hypermetrics on $7$-vertex affine bases of the Schlafli polytope.
\end{proposition}
Note that the above computation proves again that the list of $14$ orbits of hypermetric facets is complete. If not, there would exist an hypermetric facet that is violated by one extreme ray belonging to the $29$ found orbits, but this would imply that the Schlafli polytope or the $1$-simplex have interior lattice points, which is false.

\begin{corollary}
The only extreme Delaunay polytopes of dimension at most six are the $1$-simplex and the Schlafli polytope.
\end{corollary}
This method computes precisely the difference between $HYP_7$ and $CUT_7$.\\
The observed correspondence between the $37107$ non-hypermetric facets of $CUT_7$ and the $37107$ non-cut extreme rays of $HYP_7$ is presented in Table \ref{tab:tabl1}.
\begin{table}
\scriptsize
\begin{center}
\begin{tabular}{|c|ccccccccccccccccccccc|}
\hline
&12&13&14&15&16&17&23&24&25&26&27&34&35&36&37&45&46&47&56&57&67\\
\hline
$O_1$&-1&-1&0&0&1&1&-1&0&1&0&1&1&0&1&0&1&-1&1&1&-1&0\\
$R_4;\overline{G}_{24}$&2&2&2&2&1&1&2&1&1&2&1&1&1&1&2&1&2&1&1&2&2\\
\hline
$O_2;\delta_{\{3, 5, 6\}}$&-1&1&0&0&-1&1&1&0&-1&0&1&-1&0&1&0&-1&1&1&1&1&0\\
$R_5;\overline{G}_{4}$&2&1&2&1&2&1&1&1&2&1&1&2&1&1&1&2&1&1&1&1&1\\
\hline
$O_3;\delta_{\{3, 5, 4\}}$&-1&1&0&0&1&1&1&0&-1&0&1&1&0&-1&0&1&1&-1&-1&1&0\\
$R_6;\overline{G}_{23}$&2&1&1&1&1&1&1&2&2&2&1&1&1&2&1&1&1&2&2&1&2\\
\hline
\hline
$O_4$&-1&-1&-1&1&1&1&-1&0&0&1&1&0&1&0&1&1&1&1&0&-1&-1\\
$R_7;\overline{G}_{25}$&2&2&2&1&1&2&2&1&2&1&1&1&1&2&1&1&1&1&1&2&2\\
\hline
$O_5;\delta_{\{3,7\}}$&-1&1&-1&1&1&-1&1&0&0&1&-1&0&-1&0&1&1&1&-1&0&1&1\\
$R_8;\overline{G}_{5}$&2&1&2&1&1&1&1&1&2&1&2&2&2&1&1&1&1&2&1&1&1\\
\hline
$O_6;\delta_{\{2,3,7\}}$&1&1&-1&1&1&-1&-1&0&0&-1&1&0&-1&0&1&1&1&-1&0&1&1\\
$R_9;\overline{G}_{26}$&1&1&2&1&1&1&2&2&1&2&1&2&2&1&1&1&1&2&1&1&1\\
\hline
$O_7;\delta_{\{1,5,6\}}$&1&1&1&1&1&-1&-1&0&0&-1&1&0&-1&0&1&-1&-1&1&0&1&1\\
$R_{10};\overline{G}_{1}$&1&1&1&1&1&1&2&1&1&2&1&1&2&1&1&2&2&1&1&1&1\\
\hline
\hline
$O_8$&-1&-1&-1&0&1&2&-1&0&1&1&2&0&1&1&2&-1&1&1&0&-1&-2\\
$R_{11};\overline{G}_{22}$&2&2&2&2&2&1&1&1&1&1&1&1&1&1&1&2&1&1&1&2&2\\
\hline
$O_9;\delta_{\{1,4,6\}}$&1&1&-1&0&1&-2&-1&0&1&-1&2&0&1&-1&2&1&1&-1&0&-1&2\\
$R_{12};\overline{G}_{21}$&1&1&2&1&2&2&1&2&1&2&1&2&1&2&1&1&1&2&2&2&1\\
\hline
$O_{10};\delta_{\{5\}}$&-1&-1&-1&0&1&2&-1&0&-1&1&2&0&-1&1&2&1&1&1&0&1&-2\\
$R_{13};\overline{G}_{20}$&2&2&2&1&2&1&1&1&2&1&1&1&2&1&1&1&1&1&2&1&2\\
\hline
$O_{11};\delta_{\{3,5\}}$&-1&1&-1&0&1&2&1&0&-1&1&2&0&1&-1&-2&1&1&1&0&1&-2\\
$R_{14};\overline{G}_{19}$&2&1&2&1&2&1&2&1&2&1&1&2&1&2&2&1&1&1&2&1&2\\
\hline
$O_{12};\delta_{\{1,7\}}$&1&1&1&0&-1&2&-1&0&1&1&-2&0&1&1&-2&-1&1&-1&0&1&2\\
$R_{15};\overline{G}_{7}$&1&1&1&1&1&1&1&1&1&1&2&1&1&1&2&2&1&2&1&1&1\\
\hline
$O_{13};\delta_{\{7,4,1\}}$&1&1&-1&0&-1&2&-1&0&1&1&-2&0&1&1&-2&1&-1&1&0&1&2\\
$R_{16};\overline{G}_{8}$&1&1&2&1&1&1&1&2&1&1&2&2&1&1&2&1&2&1&1&1&1\\
\hline
$O_{14};\delta_{\{6,4\}}$&-1&-1&1&0&-1&2&-1&0&1&-1&2&0&1&-1&2&1&1&-1&0&-1&2\\
$R_{17};\overline{G}_{18}$&2&2&1&2&1&1&1&2&1&2&1&2&1&2&1&1&1&2&2&2&1\\
\hline
\hline
$O_{15}$&-1&-1&-2&1&1&2&0&-1&1&1&2&-2&1&1&1&2&2&3&-1&-2&-2\\
$R_{18};\overline{G}_{14}$&2&2&2&2&2&1&1&1&1&1&1&2&1&1&2&1&1&1&1&2&2\\
\hline
$O_{16};\delta_{\{5,3\}}$&-1&1&-2&-1&1&2&0&-1&-1&1&2&2&1&-1&-1&-2&2&3&1&2&-2\\
$R_{19};\overline{G}_{15}$&2&1&2&1&2&1&2&1&2&1&1&1&1&2&1&2&1&1&2&1&2\\
\hline
$O_{17};\delta_{\{5,4\}}$&-1&-1&2&-1&1&2&0&1&-1&1&2&2&-1&1&1&2&-2&-3&1&2&-2\\
$R_{20};\overline{G}_{17}$&2&2&1&1&2&1&1&2&2&1&1&1&2&1&2&1&2&2&2&1&2\\
\hline
$O_{18};\delta_{\{7,2,6\}}$&-1&1&2&-1&-1&2&0&1&-1&-1&2&-2&1&1&-1&2&2&-3&-1&2&2\\
$R_{21};\overline{G}_{13}$&2&1&1&1&1&1&2&2&2&2&1&2&1&1&1&1&1&2&1&1&1\\
\hline
$O_{19};\delta_{\{7,4,1\}}$&1&1&-2&-1&-1&2&0&1&1&1&-2&2&1&1&-1&-2&-2&3&-1&2&2\\
$R_{22};\overline{G}_{6}$&1&1&2&1&1&1&1&2&1&1&2&1&1&1&1&2&2&1&1&1&1\\
\hline
$O_{20};\delta_{\{1,7\}}$&1&1&2&-1&-1&2&0&-1&1&1&-2&-2&1&1&-1&2&2&-3&-1&2&2\\
$R_{23};\overline{G}_{2}$&1&1&1&1&1&1&1&1&1&1&2&2&1&1&1&1&1&2&1&1&1\\
\hline
$O_{21};\delta_{\{4,5,6\}}$&-1&-1&2&-1&-1&2&0&1&-1&-1&2&2&-1&-1&1&2&2&-3&-1&2&2\\
$R_{24};\overline{G}_{16}$&2&2&1&1&1&1&1&2&2&2&1&1&2&2&2&1&1&2&1&1&1\\
\hline
\hline
$O_{22}$&-1&-1&-2&1&2&3&-1&-2&1&2&3&-2&1&2&3&2&3&5&-2&-3&-5\\
$R_{25};\overline{G}_{11}$&1&1&2&2&1&1&1&2&2&1&1&2&2&1&1&1&2&1&2&2&2\\
\hline
$O_{23};\delta_{\{3,6\}}$&-1&1&-2&1&-2&3&1&-2&1&-2&3&2&-1&2&-3&2&-3&5&2&-3&5\\
$R_{26};\overline{G}_{10}$&1&2&2&2&2&1&2&2&2&2&1&1&1&1&2&1&1&1&1&2&1\\
\hline
$O_{24};\delta_{\{7,4\}}$&-1&-1&2&1&2&-3&-1&2&1&2&-3&2&1&2&-3&-2&-3&5&-2&3&5\\
$R_{27};\overline{G}_{9}$&1&1&1&2&1&2&1&1&2&1&2&1&2&1&2&2&1&1&2&1&1\\
\hline
$O_{25};\delta_{\{5\}}$&-1&-1&-2&-1&2&3&-1&-2&-1&2&3&-2&-1&2&3&-2&3&5&2&3&-5\\
$R_{28};\overline{G}_{12}$&1&1&2&1&1&1&1&2&1&1&1&2&1&1&1&2&2&1&1&1&2\\
\hline
$O_{26};\delta_{\{5,4\}}$&-1&-1&2&-1&2&3&-1&2&-1&2&3&2&-1&2&3&2&-3&-5&2&3&-5\\
$R_{29};\overline{G}_{3}$&1&1&1&1&1&1&1&1&1&1&1&1&1&1&1&1&1&2&1&1&2\\
\hline
\end{tabular}
\caption{Non-hypermetric facets of $CUT_7$ and non-cut extreme rays of $HYP_7$}\label{tab:tabl1}
\end{center}
\end{table}
The first line of Table \ref{tab:tabl1} indicates $ij$ position of the vector, defining facets and generators of extreme rays. By double line we separate $26$ pairs (facet and corresponding extreme ray) into five {\em switching} classes. Two facets $F$ and $F'$ of $CUT_7$ are called {\em switching equivalent} if there exist
\begin{equation*}
\begin{array}{c}
S\subset \{1, \dots, 7\}\mbox{,~such~that~}F(\delta_S)=0,\\
F_{ij}=-F'_{ij}\mbox{~if~}\vert S\cap \{i,j\}\vert =1\mbox{~and~} F_{ij}=F'_{ij} \mbox{,~otherwise}.
\end{array}
\end{equation*}
See section $9$ of \cite{DGL92} for details on the switchings in this case. In the first column of Table \ref{tab:tabl1} is given, for each of five switching classes, the cut $\delta_S$ such that corresponding facet is obtained by the switching by $\delta_S$ from the first facet of the class. The non-hypermetric orbits of facets of $CUT_7$ are indicated by $O_i$ and the corresponding non-cut orbits of extreme rays of $HYP_7$ are indicated by $R_{i+3}$. For any extreme ray we indicate also the corresponding graph $\overline{G}_j$ (in terms of \cite{DGL92} and chapter $16$ of \cite{DL}).\\
The five switching classes of Table \ref{tab:tabl1} correspond, respectively, to the following five classes of non-hypermetric facets of $CUT_7$, in terms of \cite{DDL} and chapter $30$ of \cite{DL}: {\em parachute} facets $P1-P3$; {\em cycle} facets $C1,C4-C6$;  {\em Grishukhin} facets $G1-G7$; cycle facets $C2, C7-C12$; cycle facets $C3,C13-C16$.\\
\cite{DG} consider extreme rays of $HYP_n$ which corresponds, moreover, to the path-metric of a graph; the Delaunay polytope, generated by such hypermetrics belongs to an integer lattice and, moreover, to a root lattice. They found, amongst $26$ non-cut orbits of extreme rays of $HYP_7$, exactly twelve which are graphic: $R_{4}$, $R_{5}$, $R_{8}$, $R_{9}$, $R_{10}$, $R_{15}$, $R_{16}$, $R_{17}$, $R_{22}$, $R_{23}$, $R_{24}$, $R_{29}$. For example, $R_{10}$, $R_{23}$ and $R_{29}$ correspond to graphic hypermetrics on $K_7-C_5$, $K_7-P_4$ and $K_7-P_3$, respectively. Three of above twelve extreme hypermetrics correspond to polytopal graphs: $3$-polytopal graph corresponding to $R_4$ and $4$-polytopal graphs $K_7-C_5$, $K_7-P_4$. Remark also that the footnote and figures on pp.~242-243 of \cite{DL} mistakingly attribute the graph $\overline{G}_{18}$ to the class $q=11$ (fourth class in our terms); in fact, it belongs to the class $q=12$ (our third class) as it was rightly given originally in \cite{DGL92}.

\section{Adjacency properties of the skeleton of $HYP_7$ and of its dual}\label{geometry-of-hypermetric-cone}
\noindent We start with Table \ref{tab:tabl2} giving incidence between extreme rays and orbits of facets, i.e. on the place $ij$ is the number of facets from the orbit $F_j$ containing fixed extreme ray from the orbit $R_i$.\\
\begin{table}
\scriptsize
\begin{center}
\begin{tabular}{|c|cccccccccccccc|c|c|}
\hline
&$F_{1}$&$F_{2}$&$F_{3}$&$F_{4}$&$F_{5}$&$F_{6}$&$F_{7}$&$F_{8}$&$F_{9}$&$F_{10}$&$F_{11}$&$F_{12}$&$F_{13}$&$F_{14}$&Inc.&$|R_i|$\\
\hline
$R_{1}$&90&150&150&180&20&15&15&180&30&30&180&180&120&120&1460&7\\
$R_{2}$&80&130&80&220&20&60&0&180&40&10&240&100&320&10&1490&21\\
$R_{3}$&75&126&96&180&18&36&12&156&30&12&162&132&240&84&1359&35\\
$R_{4}$&13&7&0&0&0&0&0&0&0&0&0&0&0&0&20&2520\\
$R_{5}$&14&6&0&0&0&0&0&0&0&0&0&0&0&0&20&2520\\
$R_{6}$&13&7&0&0&0&0&0&0&0&0&0&0&0&0&20&2520\\
$R_{7}$&14&5&0&0&1&0&0&0&0&0&0&0&0&0&20&2520\\
$R_{8}$&15&4&0&0&1&0&0&0&0&0&0&0&0&0&20&1260\\
$R_{9}$&14&5&0&0&1&0&0&0&0&0&0&0&0&0&20&1260\\
$R_{10}$&15&5&0&0&0&0&0&0&0&0&0&0&0&0&20&252\\
$R_{11}$&11&7&1&1&0&0&0&0&0&0&0&0&0&0&20&2520\\
$R_{12}$&11&7&0&2&0&0&0&0&0&0&0&0&0&0&20&2520\\
$R_{13}$&12&6&2&0&0&0&0&0&0&0&0&0&0&0&20&2520\\
$R_{14}$&11&7&0&2&0&0&0&0&0&0&0&0&0&0&20&2520\\
$R_{15}$&12&7&0&1&0&0&0&0&0&0&0&0&0&0&20&1260\\
$R_{16}$&12&6&0&2&0&0&0&0&0&0&0&0&0&0&20&1260\\
$R_{17}$&12&6&2&0&0&0&0&0&0&0&0&0&0&0&20&630\\
$R_{18}$&10&6&0&2&1&0&0&1&0&0&0&0&0&0&20&2520\\
$R_{19}$&11&5&1&1&1&0&0&1&0&0&0&0&0&0&20&2520\\
$R_{20}$&10&6&0&2&1&1&0&0&0&0&0&0&0&0&20&1260\\
$R_{21}$&10&6&1&1&1&1&0&0&0&0&0&0&0&0&20&840\\
$R_{22}$&11&6&1&0&1&0&0&1&0&0&0&0&0&0&20&840\\
$R_{23}$&11&6&0&2&1&0&0&0&0&0&0&0&0&0&20&420\\
$R_{24}$&11&6&2&0&0&0&1&0&0&0&0&0&0&0&20&420\\
$R_{25}$&7&6&1&3&0&1&0&1&0&0&0&0&1&0&20&840\\
$R_{26}$&8&5&2&2&0&0&0&2&0&0&1&0&0&0&20&630\\
$R_{27}$&8&6&0&3&0&0&0&2&0&0&0&1&0&0&20&420\\
$R_{28}$&8&6&4&0&0&0&1&0&0&0&0&0&0&1&20&210\\
$R_{29}$&8&6&0&4&0&2&0&0&0&0&0&0&0&0&20&105\\
\hline
$|F_i|$&105&210&210&420&35&105&21&420&105&42&630&420&840&210&3773&37170\\
\hline
\end{tabular}
\caption{Incidence between extreme rays and orbits of facets of $HYP_7$}\label{tab:tabl2}
\end{center}
\end{table}
It turns out, curiously, that each of $19$-dimensional hypermetric faces $F_{11}$-$F_{14}$ of $21$-dimensional cone $CUT_7$ (which became simplex facets in $HYP_7$) is the intersection of a triangle facet and one of cycle facets 
corresponding, respectively, to orbits $O_{23}$, $O_{24}$, $O_{22}$, $O_{25}$ of Table \ref{tab:tabl1}.

The {\it skeleton graph} of $HYP_7$ is the graph 
whose nodes are the extreme rays of $HYP_7$ and whose edges
are the pairs of adjacent extreme rays.
The {\it ridge graph} of $HYP_7$ is the graph with node set being
the set of facets of $HYP_7$ and with an edge between two facets if
they are adjacent on $HYP_7$.

\begin{proposition}
We have the following properties of adjacency of extreme rays:
\begin{enumerate}
\item[(i)] The restriction of the skeleton of $HYP_7$ on the union of cut orbits $R_1\cup R_2\cup R_3$ is the complete graph.
\item[(ii)] Every non-cut extreme ray of $HYP_7$ has adjacency $20$ (namely, it is adjacent to $20$ cuts lying on corresponding non-hypermetric facet of $CUT_7$); see on Table \ref{tab:tabl3} the distribution of those $20$ cuts amongst the cut orbits.
\item[(iii)] Any two simplex extreme rays are non-adjacent; any simplex extreme ray (i.e. non-cut ray) has local graph (i.e. the restriction of the skeleton on the set of its neighbors) $K_{20}$.
\item[(iv)] The diameter of the skeleton graph of $HYP_7$ is $3$.
\end{enumerate}
\end{proposition}
\begin{table}
\begin{center}
\scriptsize
\begin{tabular}{|c|ccc|c|c|}
\hline
&$R_1$&$R_2$&$R_3$&Adj.&$|R_i|$\\
\hline
$R_{1}$& 6& 21& 35& 15662&7\\
$R_2$& 7& 20& 35& 12532&21\\
$R_3$& 7& 21& 34& 10664&35\\
\hline
$R_4$& 3& 6&  11& 20&2520\\
$R_5$& 4& 7&   9& 20&2520\\
$R_6$& 3& 7& 10& 20&2520\\
\hline
$R_7$& 3& 7& 10& 20&2520\\
$R_8$& 4& 7& 9& 20&1260\\
$R_9$& 3& 8& 9& 20&1260\\
$R_{10}$& 5& 5& 10& 20&252\\
\hline
$R_{11}$& 3& 6& 11& 20&2520\\
$R_{12}$& 2& 8& 10& 20&2520\\
$R_{13}$& 4& 5& 11& 20&2520\\
$R_{14}$& 2& 8& 10& 20&2520\\
$R_{15}$& 4& 7& 9& 20&1260\\
$R_{16}$& 3& 9& 8& 20&1260\\
$R_{17}$& 4& 4& 12& 20&630\\
\hline
$R_{18}$& 2& 8& 10& 20&2520\\
$R_{19}$& 3& 7& 10& 20&2520\\
$R_{20}$& 1& 10& 9& 20&1260\\
$R_{21}$& 2& 9& 9& 20&840\\
$R_{22}$& 4& 6& 10& 20&840\\
$R_{23}$& 4& 7& 9& 20&420\\
$R_{24}$& 5& 1& 14& 20&420\\
\hline
$R_{25}$& 1& 9& 10& 20&840\\
$R_{26}$& 2& 8& 10& 20&630\\
$R_{27}$& 3& 6& 11& 20&420\\
$R_{28}$& 5& 1& 14& 20&210\\
$R_{29}$& 2& 12& 6& 20&105\\
\hline
$|R_j|$&7&21&35&&37170\\
\hline
\end{tabular}
\caption{Adjacencies between extreme rays and orbits of extreme rays of $HYP_7$}\label{tab:tabl3}
\end{center}
\end{table}

\begin{proposition}
We have the following properties of adjacency of facets of $HYP_7$:
\begin{enumerate}
\item[(i)] See Table \ref{tab:tabl4}, where on the place $ij$ we have the number of facets from orbit $F_j$ adjacent to fixed facet of orbit $F_i$.
\item[(ii)] Any two simplex facets are non-adjacent; any simplex facet (i.e. from $F_9$--$F_{14}$) have local graph $K_{20}$.
\item[(iii)] The diameter of the ridge graph of $HYP_7$ is $3$.
\item[(iv)] The symmetry group of the skeleton of dual $HYP_7$ is the symmetric group $Sym(7)$.
\end{enumerate}
\end{proposition}

\begin{table}
\begin{center}
\scriptsize
\begin{tabular}{|c|cccccccccccccc|c|c|}
\hline
&$F_{1}$&$F_{2}$&$F_{3}$&$F_{4}$&$F_{5}$&$F_{6}$&$F_{7}$&$F_{8}$&$F_{9}$&$F_{10}$&$F_{11}$&$F_{12}$&$F_{13}$&$F_{14}$&Adj.&$|F_i|$\\
\hline
$F_{1}$& 86& 168& 110& 216& 35& 56& 13& 196& 14& 6& 54& 36& 64& 18&1072&105\\
$F_{2}$& 84& 116& 62& 114& 3& 5& 1& 18& 0& 0& 15& 12& 24& 6&460&210\\
$F_{3}$& 55& 62& 9& 20& 1& 1& 1& 4& 1& 1& 6& 0& 4& 4&169&210\\
$F_{4}$& 54& 57& 10& 25& 2& 2& 0& 6& 1& 0& 3& 3& 6& 0&169&420\\
$F_{5}$& 105& 18& 6& 24& 0& 3& 0& 12& 0& 0& 0& 0& 0& 0&168&35\\
$F_{6}$& 56& 10& 2& 8& 1& 2& 0& 8& 0& 0& 0& 0& 8& 0&95&105\\
$F_{7}$& 65& 10& 10& 0& 0& 0& 0& 0& 0& 0& 0& 0& 0& 10&95&21\\
$F_{8}$& 49& 9& 2& 6& 1& 2& 0& 5& 0& 0& 3& 2& 2& 0&81&420\\
$F_{9}$& 14& 0& 2& 4& 0& 0& 0& 0& 0& 0& 0& 0& 0& 0&20&105\\
$F_{10}$& 15& 0& 5& 0& 0& 0& 0& 0& 0& 0& 0& 0& 0& 0&20&42\\
$F_{11}$& 9& 5& 2& 2& 0& 0& 0& 2& 0& 0& 0& 0& 0& 0&20&630\\
$F_{12}$& 9& 6& 0& 3& 0& 0& 0& 2& 0& 0& 0& 0& 0& 0&20&420\\
$F_{13}$& 8& 6& 1& 3& 0& 1& 0& 1& 0& 0& 0& 0& 0& 0&20&840\\
$F_{14}$& 9& 6& 4& 0& 0& 0& 1& 0& 0& 0& 0& 0& 0& 0&20&210\\
\hline
$|F_j|$&105&210&210&420&35&105&21&420&105&42&630&420&840&210&&3773\\
\hline
\end{tabular}
\caption{Adjacency between facets and orbits of facets of $HYP_7$}\label{tab:tabl4}
\end{center}
\end{table}

\section{Final remarks}\label{further-research}
\noindent In order to find extreme rays of $HYP_8$ the same methods will, probably, work with more computational difficulties but in dimension $n\geq 9$ polyhedral methods may fail.\\
The list of $374$ orbits of non-cut extreme rays of $HYP_8$ (containing $7126560$ extreme rays), found in \cite{Du}, will be confronted there with the list of at least $2169$ orbits of facets of $CUT_8$, found in \cite{CR}. Exactly $55$ of above $374$ orbits corresponds to path-metrics of a graph. It was shown in \cite{DG} that any graph, such that its path-metric lies on an extreme ray of a $HYP_n$, is a subgraph of the skeleton of Gosset or Schafli polytopes.\\
It turns out (it can also be found in chapter $28$ of \cite{DL}) that exactly $26$ of those orbits consist of hypermetric inequalities; ten are $0$-extensions of the hypermetric facets of $CUT_7$ and $16$ come from the following vectors $b$ (see (\ref{Definition-of-hypermetrics})):
\begin{equation*}
\begin{array}{ccc}
(2, 1, 1, 1, -1, -1, -1, -1),&	(3, 1, 1, 1, -1, -1, -1, -2),& (2, 2, 1, 1, -1, -1, -1, -2),\\
(4, 1, 1, -1, -1, -1, -1, -1),&	(3, 2, 2, -1, -1, -1, -1, -2),
\end{array}
\end{equation*}
representing, respectively, switching classes of sizes $2$, $4$, $3$, $2$, $5$.\\
There is one to one correspondance between non-hypermetric facets of $CUT_7$ and non-cut extreme rays of $HYP_7$; namely every such facet is violated by exactly one such ray.
Apropos, there is one to one correspondance between the ten non-cut extreme rays (in fact, path metric $K_{2\times 3}$) of $MET_5$ and the ten non-triangle facets (in fact, pentagonal) of $CUT_5$; namely every non-cut extreme ray is violated by a non-triangle facet. There is no such thing between $MET_n$ and $CUT_n$ for $n>5$ but we hope that the correspondance exist for $CUT_8$ and $HYP_8$.\\
Another direction for further study is to find {\em all} faces of $HYP_7$. While the extreme rays of $HYP_n$ yields the extreme Delaunay polytopes of dimension $n-1$, the study of all faces of $HYP_n$ will provide the list of all (combinatorial types of) Delaunay polytopes of dimension less or equal $n-1$. See \cite{BK} for description of the method and computations for $n\leq 4$; in fact, Kononenko (submitted) found all Delaunay polytopes of dimension five.


\begin{thebibliography}{99}

\bibitem[As82]{As}
P. Assouad, {\em Sous-espaces de $L^1$ et in\'egalit\'es hyperm\'etriques}, Compte Rendus de l'Acad\'emie des Sciences de Paris, {\bf 294(A)} (1982) 439--442.


\bibitem[B70]{B70}
E.P. Baranovski, {\em Simplexes of $L$-subdivisions of euclidean spaces}, Mathematical Notes, {\bf 10} (1971) 827-834.

\bibitem[Ba99]{Ba}
E.P. Baranovskii, {\em The conditions for a simplex of $6$-dimensional lattice to be $L$-simplex}, (in russian) Nauchnyie Trudi Ivanovo state university. Mathematica {\bf 2} (1999) 18--24.

\bibitem[BK00]{BK}
E.P. Baranovskii and P.G. Kononenko, {\em A method of deducing $L$-polyhedra for $n$-lattices}, Mathematical Notes {\bf 68-6} (2000) 704--712.

\bibitem[ChRe98]{CR}
T. Christof and G. Reinelt, {\em Decomposition and Parallelization Techniques for Enumerating the Facets of 0/1-Polytopes}, Preprint, Univ Heidelberg, 1998.

\bibitem[CS]{CS}
J.H. Conway and N.J.A. Sloane, {\em Sphere Packings, Lattices and Groups}, volume 290 of {\em Grundlehren der mathematischen Wissenschaften}, Springer Verlag.

\bibitem[DDL94]{DDL}
C. De Simone, M. Deza, and M. Laurent, {\em Collapsing and lifting for the cut cone}, Discrete Mathematics, {\bf 127} (1994), 105--140.

\bibitem[DeDe94]{DD1}
A. Deza and M. Deza, {\em The ridge graph of the metric polytope and some 
relatives}, in T.Bisztriczky, P.McMullen, R.Schneider and A.Ivic Weiss eds
{\em Polytopes: Abstract, Convex and Computational} (1994) 359--372.

\bibitem[DG93]{DG}
M. Deza and V.P. Grishukhin, {\em Hypermetric graphs}, The Quarterly Journal of Mathematics Oxford, (2) {\bf 44} (1993) 399--433.

\bibitem[DGL92]{DGL92}
M. Deza, V.P. Grishukhin, and M. Laurent, {\em Extreme hypermetrics and L-polytopes}, in G.Hal\'asz et al. eds {\em Sets, Graphs and Numbers, Budapest (Hungary), 1991}, {\bf 60} {\em Colloquia Mathematica Societatis J\'anos Bolyai}, (1992) 157--209.

\bibitem[DGL93]{DGL93}
M.Deza, V.P. Grishukhin, and M. Laurent, {\em The hypermetric cone is polyhedral}, Combinatorica, {\bf 13} (1993) 397--411.

\bibitem[DeLa97]{DL}
M. Deza and M. Laurent, {\em Geometry of cuts and metrics}, Springer--Verlag,
 Berlin, 1997.

\bibitem[Du]{Du}
M. Dutour, {\em The Gosset polytope and the hypermetric cone on eight vertices}, in preparation.

\bibitem[Fu]{Fu}
K. Fukuda, {\em http://www.ifor.math.ethz.ch/\texttt{\~}fukuda/cdd\_home/cdd.html}

\bibitem[Gr90]{G90}
V.P. Grishukhin, {\em All facets of the cut cone $C_n$ for $n=7$ are known}, European Journal of Combinatorics, {\bf 11} (1990) 115--117.

\bibitem[Gr92]{GrMet}
V.P. Grishukhin, {\em Computing extreme rays of the metric cone for seven points}, European Journal of Combinatorics, {\bf 13} (1992) 153--165.


\bibitem[RB]{BR}
S.S. Ryshkov and E.P. Baranovskii, {\em The Repartitioning Complexes in n-dimensional Lattices (with Full Description for $n\leq 6$)}, in Voronoi's impact on modern science, Book 2, Institute of Mathematics, Kyiv (1998) 115-124.

\bibitem[MK]{MK}
B. McKay, {\em http://cs.anu.edu.au/people/bdm/nauty/}

\bibitem[Vo]{Vo}
G.F. Voronoi, {\em Nouvelles applications des param\`etres continus \`a la th\'eorie des formes quadratiques - Deuxi\`eme m\'emoire}, J. f\"ur die reine und angewandte Mathematik, {\bf 134} (1908) 198-287, {\bf 136} (1909) 67-178.


\end{thebibliography}
\end{document}